\documentclass[12pt]{article}
\usepackage{amssymb}
\catcode`\@=11
\@addtoreset{equation}{section}

\catcode`\@=12

\newtheorem{Theorem}{Theorem}[section]
\newtheorem{Definition}{Definition}[section]

\newtheorem{Corollary}{Corollary}[section]
\newtheorem{Remark}{Remark}[section]

\def\IR{{\mathbb R}}

\title{Completely linear degeneracy for quasilinear hyperbolic systems in several space variables}
\author{De-Xing Kong $\quad$ and $\quad$ Chang-Hua Wei\footnote{Corresponding
author: wch\_19861125@163.com.} \\
Department of Mathematics, Zhejiang University\\ Hangzhou 310027,
China}
\date{}
\begin{document}
\maketitle
\begin{abstract} In this paper, we introduce the concept of completely linear degeneracy for quasilinear hyperbolic systems in several space variables, and then get an interesting property for multidimensional hyperbolic conservation laws. Some examples and applications are given at last.

\vskip 3mm

\end{abstract}
\section{Introduction}
For quasilinear hyperbolic systems in several space variables
\begin{equation}
\frac {\partial u}{\partial t}+\sum^m_{j=1}A_j(u)\frac {\partial u}{\partial x_j}=0,
\end{equation}
where $u=(u_{1},\cdots,u_{n})^{T}$ is the unknown vector function of $(t,x_{1},\cdots,x_{m})$ and $A_{j}(u)\,(j=1,\cdots,m)$ ia an $n\times n$ matrix with smooth elements $a_{jkl}(u)\,(k,l=1,\cdots,n)$.
The concepts of linear degeneracy and genuine nonlinearity
have been made in the following way (see \cite{ma}). This is a straightforward generalization of the case of one space dimension (see \cite{lax}). The $i$-th characteristic field of system (1.1)
is linearly degenerate, if
\begin{equation}
\nabla\lambda_i(u,\xi )\cdot r_i(u,\xi )\equiv 0,\quad \forall\; u\in\Omega,\;\;\forall\;\xi\in {\mathbb S}^{m-1};
\end{equation}
while, it is genuinely nonlinear, if
\begin{equation}
\nabla\lambda_i(u,\xi )\cdot r_i(u,\xi )\neq 0, \quad \forall\; u\in\Omega,\;\;\forall\;\xi\in {\mathbb S}^{m-1},
\end{equation}
where $\xi =(\xi_1,\cdots ,\xi_m)^T\in {\mathbb S}^{m-1}$, ${\mathbb S}^{m-1}$ is the unit ball centered at the origin in $\IR^m$, $\lambda_1\left(u,\xi\right),\cdots ,\lambda_n\left(u,\xi\right)$ are $n$ real eigenvalues of $A(u,\xi)$ and $\{r_i\left(u,\xi\right)\}^n_{i=1}$ is a complete set of right eigenvectors of $A(u,\xi)$, in which $A(u,\xi)=\sum^m_{j=1}A_j(u)\xi_j$. Here we assume that $A(u,\xi)$ has $n$ real eigenvalues. However, as pointed out in \cite{lax1}, this  generalization would for instance exclude a single equation and a system of two equations being
genuinely nonlinear (see \cite{lax1}) and thus is unsuitable.

In this paper, we investigate some basic properties enjoyed by quasilinear hyperbolic systems in several space dimensions. We firstly introduce a concept of ``completely linear degeneracy'' for quasilinear hyperbolic system in several space dimensions and then study some interesting properties enjoyed by multi-dimensional hyperbolic conservation laws. In particular, we prove that the system of hyperbolic conservation laws with two unknowns in several space dimensions must be linear if the system is completely linearly degenerate. This gives a criterion on the linearity of this kind of system. Some examples and applications are also given.

\section{Completely linear degeneracy}
We introduce
\begin{Definition}
 The $i$-th characteristic field of (1.1) is said to be completely linearly degenerate, if it holds that
\begin{equation}
\nabla_{u}\lambda_{i}(u,\xi)\cdot r_{i}(u,\xi)\equiv 0,\qquad \forall\,u\in\Omega,\quad \forall\,\xi\in\mathbb S^{m-1}
\end{equation}
and
\begin{equation}
\nabla_{u} r_{ik}(u,\xi)\cdot r_{i}(u,\xi)\equiv 0,\qquad \forall\,u\in\Omega,\quad \forall\,\xi\in\mathbb S^{m-1},\quad \forall\;k\in\{1,2,\cdots,n\},
\end{equation}
where $r_{ik}(u,\xi)$ stands for the $k$-th component of the vector $r_{i}(u,\xi)$.

We call system (1.1) completely linearly degenerate, if all characteristic fields are completely linearly degenerate.
\end{Definition}
\begin{Remark}
The quantities $\nabla_{u}\lambda_{i}(u,\xi)\cdot r_{i}(u,\xi)$ and $\nabla_{u}r_{ik}(u,\xi)\cdot r_{i}(u,\xi)\;k\in\{1,\cdots,n\}$ are invariant under any invertible smooth transformation of unknowns.
\end{Remark}
\section{An interesting property}
In this section, we consider the following hyperbolic conservation laws with two unknows:
\begin{equation}
\partial_{t}\textbf{u}+\nabla_{x}\cdot f(\textbf{u})=0,\qquad \textbf{u}=(u,v)^{T}\in\mathbb R^{2},\quad x\in\mathbb R^{n},
\end{equation}
where $\nabla_{x}=\{\partial_{x_{1}},\cdots,\partial_{x_{n}}\}$ and
$$
f=(f_{1},\cdots,f_{n}):\mathbb R^{2}\rightarrow (\mathbb R^{2})^{n}
$$
is a nonlinear smooth mapping with $f_{i}:\mathbb R^{2}\rightarrow \mathbb R^{2}$ for $i=1,\cdots,n$.

Consider plane wave solutions
$$
\textbf{u}(t,x)=\textbf{w}(t,x\cdot\xi)\qquad {\rm for}\quad \xi\in\mathbb S^{n-1}.
$$
Then $\textbf{w}(t,s)$ satisfies
$$
\partial_{t}\textbf{w}+(\nabla f(\textbf{w})\cdot\xi)\partial_{s}\textbf{w}=0,
$$
where $\nabla=\{\partial_{w_{1}},\partial_{w_{2}}\}$.

Assume that system (3.1) is hyperbolic in a state domain $\Omega$, that is to say, for every state $\textbf{u}\in\Omega$, $(\nabla f(\textbf{w})\cdot\xi)$ has two real eigenvalues $\lambda_{1}(\textbf{u},\xi),\,\lambda_{2}(\textbf{u},\xi)$
($\lambda_{1}(\textbf{u},\xi)=\lambda_{2}(\textbf{u},\xi)$ is included) and a complete set of right eigenvectors $r_{1}(\textbf{u},\xi),\,r_{2}(\textbf{u},\xi)$. For $i=1,2$ let $r_{i}(\textbf{u},\xi)=(r_{i1}(\textbf{u},\xi),\,r_{i2}(\textbf{u},\xi))^{T}$ be an eigenvector corresponding to $\lambda_{i}(\textbf{u},\xi)$, i.e.,
\begin{equation}
(\nabla f(\textbf{u})\cdot\xi)r_{i}(\textbf{u},\xi)=\lambda_{i}(\textbf{u},\xi)r_{i}(\textbf{u},\xi)
\end{equation}
and
\begin{equation}
det |r_{ij}(\textbf{u},\xi)|\neq0\qquad (i,j=1,2).
\end{equation}
With the above preparations, we have
\begin{Theorem}
The system of hyperbolic conservation laws (3.1) is completely linearly degenerate in the sense of Definition 2.1 if and only if (3.1) is a linear system.
\end{Theorem}
{\it Proof.}
The necessity is obvious, so we prove the sufficiency only. Assume that system (3.1) is completely linearly degenerate, i.e., for any state $\textbf{u}\in\Omega$ and $\xi=(\xi_{1},\cdots,\xi_{n})\in\mathbb S^{n-1}$, $\nabla f(\textbf{u},\xi)\cdot\xi$ has two real eigenvalues $\lambda_{1}(\textbf{u},\xi),\,\lambda_{2}(\textbf{u},\xi)$ and two corresponding right eigenvectors $r_{1}(\textbf{u},\xi),\,r_{2}(\textbf{u},\xi)$, which satisfy
\begin{equation}
det |r_{ij}(\textbf{u},\xi)|\neq0\qquad (i,j=1,2).
\end{equation}
Specially, if we choose $\xi=(0,\cdots,0,\displaystyle{\mathop{1}^{(i)}},0,\cdots,0)\;(i=1,\cdots,n)$, then it must hold that
$\nabla f_{i}(\textbf{u})$ has two completely linearly degenerate characteristic fields. Since
\begin{equation}
\nabla f_{i}(\textbf{u})=\left(\begin{array}{cc}
(f_{i1})_{u}& (f_{i1})_{v}\\
(f_{i2})_{u}& (f_{i2})_{v}
\end{array}\right).
\end{equation}
For convenience, denote
$$
(f_{i1})_{u}=a(u,v),\quad (f_{i1})_{v}=b(u,v),\quad (f_{i2})_{u}=c(u,v)\quad {\rm and} \quad (f_{i2})_{v}=d(u,v),
$$
then it holds that
\begin{equation}
a_{v}(u,v)=b_{u}(u,v)\quad {\rm and} \quad c_{v}(u,v)=d_{u}(u,v).
\end{equation}
By a direct calculation, we have
\begin{equation}
\lambda_{1}(u,v)=\frac{a+d+\sqrt{(a-d)^{2}+4bc}}{2}\quad {\rm and} \quad \lambda_{2}(u,v)=\frac{a+d-\sqrt{(a-d)^{2}+4bc}}{2}.
\end{equation}

Next, we prove the theorem by dividing the following four cases, suppose that in a neighborhood of some state $(u,v)\in\Omega$

Case I.\quad $(a-d)^{2}+4bc>0$, $b\neq0$ and $c\neq0$.\\
Choose the corresponding right eigenvectors as
\begin{equation}
r_{1}(u,v)=(b,\; \lambda_{1}(u,v)-a)^{T}\quad {\rm and}\quad r_{2}(u,v)=(b,\; \lambda_{2}(u,v)-a)^{T}.
\end{equation}
Then, by Definition 2.1, it holds that
\begin{equation}
\nabla \lambda_{1}(u,v)\cdot r_{1}(u,v)\equiv 0,\; \nabla b(u,v)\cdot r_{1}(u,v)\equiv0\; {\rm and} \; \nabla\left(\lambda_{1}(u,v)-a(u,v)\right)\cdot r_{1}(u,v)\equiv0
\end{equation}
and
\begin{equation}
\nabla \lambda_{2}(u,v)\cdot r_{2}(u,v)\equiv 0,\; \nabla b(u,v)\cdot r_{2}(u,v)\equiv0\;{\rm and} \; \nabla\left(\lambda_{2}(u,v)-a(u,v)\right)\cdot r_{2}(u,v)\equiv0.
\end{equation}
By (3.9) and (3.10), we have
\begin{equation}
\nabla a(u,v)\cdot r_{1}(u,v)=0\;{\rm and} \;\nabla a(u,v)\cdot r_{2}(u,v)=0
\end{equation}
and
\begin{equation}
\nabla b(u,v)\cdot r_{1}(u,v)=0\;{\rm and} \;\nabla b(u,v)\cdot r_{2}(u,v)=0.
\end{equation}
Since
$$
\det|r_{1}(u,v),r_{2}(u,v)|\neq0,\quad \forall\,(u,v)\in\Omega,
$$
by (3.11) and (3.12), it holds that
\begin{equation}
\nabla a(u,v)\equiv0\quad{\rm and} \quad \nabla b(u,v)\equiv0,
\end{equation}
i.e., $a=const.$ and $b=const.$.

On the other hand, if we choose
$$
\widetilde{r}_{1}(u,v)=(\lambda_{2}(u,v)-d(u,v),\,c(u,v))^{T},\,\widetilde{r}_{2}(u,v)=(\lambda_{2}(u,v)-d(u,v),\, c(u,v))^{T},
$$by Definition 2.1, it must hold that
\begin{equation}
\nabla c(u,v)\equiv0\quad{\rm and} \quad \nabla d(u,v)\equiv 0,
\end{equation}
i.e., $c=const.$ and $d=const.$.

Case II.\quad $(a-d)^{2}+4bc>0$, $b\equiv0$ and $c\neq0$.\\
By a direct calculation, the eigenvalues are
\begin{equation}
\lambda_{1}(u,v)=a(u,v)\quad{\rm and}\quad \lambda_{2}(u,v)=d(u,v)
 \end{equation}
 and the corresponding right eigenvectors can be chosen as
 \begin{equation}
 r_{1}(u,v)=(a(u,v)-d(u,v),\,c(u,v))^{T}\quad{\rm and} \quad r_{2}(u,v)=(0,1)^{T}.
 \end{equation}
 Then by Definition 2.1, it holds that
 \begin{equation}
 \nabla a(u,v)\cdot r_{1}(u,v)=0,\,\nabla c(u,v)\cdot r_{1}(u,v)=0,\,\nabla\left(a(u,v)-d(u,v)\right)\cdot r_{1}(u,v)=0
 \end{equation}and
 \begin{equation}
 \nabla d(u,v)\cdot r_{2}(u,v)=0.
 \end{equation}
 By (3.17) and (3.18), we have
 \begin{equation}
 \nabla d(u,v)\cdot r_{1}(u,v)=0\;{\rm and}\;\nabla d(u,v)\cdot r_{2}(u,v)=0
 \end{equation}and
 \begin{equation}
\left\{\begin{array}{c}
 a_{u}(a-d)+a_{v}c=0\\
 c_{u}(a-d)+c_{v}c=0
 \end{array}\right.
 \end{equation}
 By (3.19), we have
 \begin{equation}
 d=const.
\end{equation}
By (3.6) and (3.20), it holds that
\begin{equation}
a_{v}=0\quad{\rm and} \quad c_{v}=0.
\end{equation}
By the assumption that $a-d\neq0$, we have by (3.20)
\begin{equation}
a_{u}=0\quad{\rm and} \quad c_{u}=0.
\end{equation}
(3.22) and (3.23) imply that $a=const.$ and $c=const.$.

Case III.\quad $a(u,v)-d(u,v)\neq0$ and $b(u,v)=c(u,v)\equiv0$.\\
The eigenvalues are $$\lambda_{1}(u,v)=a(u,v)\quad{\rm and}\quad\lambda_{2}(u,v)=d(u,v)$$ respectively and the corresponding right eigenvectors can be chosen as
 $$r_{1}=(1,0)^{T}\quad{\rm and} \quad r_{2}=(0,1)^{T}.$$  Then by Definition 2.1, it holds that
\begin{equation}
a_{u}=0\quad{\rm and} \quad d_{v}=0.
\end{equation}
By (3.6), we have
\begin{equation}
a_{v}=0\quad{\rm and} \quad d_{u}=0.
\end{equation}
Thus, $a=const.$ and $d=const.$ hold.

Case IV. \quad $(a-d)^{2}+4bc=0$. i.e., $a(u,v)=d(u,v)$ and b(u,v)c(u,v)=0.\\
By a direct calculation, the eigenvalues are $$\lambda_{1}(u,v)=\lambda_{2}(u,v)=\frac{a(u,v)+d(u,v)}{2}.$$ By the hyperbolicity of system (3.1), it must hold $b(u,v)=c(u,v)\equiv0$. Then, by Case III, $a=d=const.$ holds.

By the continuity of $a,b,c$ and $d$, we can extend the neighborhood of the state $(u,v)$ to the whole state space $\Omega$. Thus, the theorem is proved.
\begin{Remark}
For the case that at some isolated states $(u,v)$, $b=0$ or $c=0$ but at other states $b\neq0$ and $c\neq0$, one can refer to Case I.
\end{Remark}
\begin{Corollary}
The solution of the Cauchy problem of multi-dimensional hyperbolic conservation laws with two unknowns exists globally if it is completely linearly degenerate in the sense of Definition 2.1.
\end{Corollary}

\section{Examples and applications}
\subsection{System of gas dynamics}
Now let us consider the system of gas dynamics in two space dimensions:
\begin{equation}\begin{array}{l}
\frac{\partial\rho}{\partial t}+\frac{\partial (\rho u)}{\partial x}+ \frac{\partial (\rho v)}{\partial y}=0,\\
\frac{\partial (\rho u)}{\partial t}+ \frac{\partial\left(\rho u^2+p\right) }{\partial x}+
\frac{\partial (\rho uv)}{\partial y}=0,\\
\frac{\partial (\rho v)}{\partial t}+\frac{\partial (\rho uv)}{\partial x}+\frac{\partial \left(\rho v^2+p\right)}{\partial y}=0,\\
\frac{\partial \left(\rho/2\left(u^2+v^2\right)+\rho e\right)}{\partial t}
+\frac{\partial \left[\left(\rho /2 \left( u^2+v^2\right)+\rho e+p\right)u \right] }{\partial x}+\frac{\partial \left[ \left( \rho /2 \left( u^2+v^2\right)+\rho e+p\right)v \right] }{\partial y}=0.  \end{array}\end{equation}
For classical solutions, it is equivalent to
\begin{equation}\begin{array}{l}
\frac{\partial \rho}{\partial t}+\frac{\partial (\rho u)}{\partial x}+ \frac{\partial (\rho v)}{\partial y}=0,\\
\frac{\partial u}{\partial t}+ u\frac{\partial u}{\partial x}+v\frac{\partial u}{\partial y}+\frac{1}{\rho}\frac{\partial p}{\partial x}=0,\\
\frac{\partial v}{\partial t}+u\frac{\partial v}{\partial x}+v\frac{\partial v}{\partial y}+\frac{1}{\rho}\frac{\partial p}{\partial y}=0,\\
\frac{\partial s}{\partial t}+u\frac{\partial s}{\partial x}+v\frac{\partial s}{\partial y}=0.
\end{array}\end{equation}
We can easily calculate
$$\lambda_1=u\xi_1 + v\xi_2 -\sqrt{p_{\rho}}\sqrt{\xi^2_1+\xi^2_2},$$
$$\lambda_2=\lambda_3=u\xi_1+v\xi_2,$$
$$\lambda_4=u\xi_1+v\xi_2+\sqrt{p_{\rho}}\sqrt{\xi^2_1+\xi^2_2},$$
$$\nabla\lambda_1 r_1 = \nabla\lambda_4 r_4 \neq 0,$$
$$\nabla\lambda_2 r_2 = \nabla\lambda_3 r_3 \equiv 0.$$
Thus, $\lambda_1$ and $\lambda_4$ are genuinely nonlinear. For classical solutions, the system is also equivalent to
\begin{equation}\begin{array}{l}
\frac{\partial \rho}{\partial t}+\frac{\partial (\rho u)}{\partial x}+ \frac{\partial (\rho v)}{\partial y}=0,\\
\frac{\partial u}{\partial t}+ u\frac{\partial u}{\partial x}+v\frac{\partial u}{\partial y}+\frac{1}{\rho}\frac{\partial p}{\partial x}=0,\\
\frac{\partial v}{\partial t}+ u\frac{\partial v}{\partial x}+v\frac{\partial v}{\partial y}+\frac{1}{\rho}\frac{\partial p}{\partial y}=0,\\
\frac{\partial p}{\partial t}+ u\frac{\partial p}{\partial x}+v\frac{\partial p}{\partial y}+\frac{\partial p}{\partial \rho}\left( \frac{\partial u}{\partial x}+\frac{\partial v}{\partial y}\right)=0.\end{array}\end{equation}
In this case
$$r_2=(1,0,0,0),\quad r_3=(0, -\xi_2, \xi_1, 0),$$
so $\lambda_2 =\lambda_3$ are completely linearly degenerate in the sense of Definition 2.1.
\subsection{Relativistic torus in the Minkowski space $\mathbb R^{1+n}$}
The motion of relativistic torus in $\mathbb R^{1+n}$ corresponds to the three-dimensional time-like extremal submanifolds in $\mathbb R^{1+n}$ and plays an important role in both mathematics and physics, which is studied by Huang and Kong \cite{HK}. The equations governing the motion of relativistic torus in the Minkowski space $\mathbb R^{1+n}$ are
\begin{equation}
\begin{array}{ll}
(|x_{\alpha}|^{2}|x_{\beta}|^{2})x_{tt}+2(<x_{\alpha},x_{\beta}><x_{t},x_{\beta}>-<x_{t},x_{\alpha}>|x_{\beta}|^{2})
x_{t\alpha}\\
\quad+\left((|x_{t}|^{2}-1)|x_{\alpha}|^{2}-<x_{t},x_{\alpha}>^{2}\right)x_{\beta\beta}
+2\left(<x_{t},x_{\alpha}><x_{\alpha},x_{\beta}>-<x_{t},x_{\beta}>|x_{\alpha}|^{2}\right)x_{t\beta}\\
\quad+\left((|x_{t}|^{2}-1)|x_{\beta}|^{2}-<x_{t},x_{\beta}>^{2}\right)x_{\alpha\alpha}\\
\quad+
2\left(<x_{t},x_{\alpha}><x_{t},x_{\beta}>-(|x_{t}|^{2}-1)<x_{\alpha},x_{\beta}>\right)x_{\alpha\beta}=0,
\end{array}
\end{equation}
where $t,\alpha,\beta$ are parameter coordinates and $x=x(t,\alpha,\beta)=\left(x_{1}(t,\alpha,\beta),\cdots,x_{n}(t,\alpha,\beta)\right)$. $<,>$ denotes the inner product in Euclidean space.

Let
\begin{equation}
u=x_{t},\quad v=x_{\alpha},\quad w=x_{\beta}
\end{equation}
and
\begin{equation}
U=(u^{T},v^{T},w^{T})^{T},
\end{equation}
where $u^{T}=(u_{1},\cdots,u_{n})$, $v^{T}=(v_{1},\cdots,v_{n})$ and $w^{T}=(w_{1},\cdots,w_{n})$. Then equation (4.4) can be equivalently rewritten as
\begin{equation}
A_{0}(U)\frac{\partial U}{\partial t}+A_{1}(U)\frac{\partial U}{\partial \alpha}+B_{1}(U)\frac{\partial U}{\partial \beta}=0
\end{equation}
for classical solutions, where
\begin{equation}
A_{0}(U)=\left[\begin{array}{ccc}
-aI&0&0\\
0& I&0\\
0&0& I
\end{array}\right],\quad
A_{1}(U)=\left[\begin{array}{ccc}
-2bI&-eI&0\\
-I&0&0\\
0&0&0
\end{array}\right],
\quad
B_{1}(U)=\left[
\begin{array}{ccc}
-2cI&-2dI&-p I\\
0&0&0\\
-I&0&0
\end{array}\right],
\end{equation}
in which
\begin{equation}
\begin{array}{cc}
a=<v,w>^{2}-|v|^{2}|w|^{2},\quad b=<u,v>|w|^{2}-<u,w><v,w>,\\
 d=<v,w>(|u|^{2}-1)-<u,v><u,w>,\quad c=<u,w>|v|^{2}-<u,v><v,w>,\\
 e=<u,w>^{2}-|w|^{2}(|u|^{2}-1),\quad p=<u,v>^{2}-|v|^{2}(|u|^{2}-1)
 \end{array}
\end{equation}
and $I$ represents an $n\times n$ unit matrix. Since the submanifold is timelike, we have $a<0$ and $A_{0}(U)$ is nonsingular, then we can rewrite equation (4.7) as
\begin{equation}
\frac{\partial U}{\partial t}+A(U)\frac{\partial U}{\partial \alpha}+B(U)\frac{\partial U}{\partial \beta}=0,
\end{equation}
where
\begin{equation}
A(U)=\left[\begin{array}{ccc}
\frac{2b}{a}I&\frac{e}{a}I&0\\
-I&0&0\\
0&0&0
\end{array}\right],\quad
B(U)=\left[
\begin{array}{ccc}
\frac{2c}{a}I&\frac{2d}{a}I&\frac{p}{a} I\\
0&0&0\\
-I&0&0
\end{array}\right].
\end{equation}

Set
\begin{equation}
M(U)=\xi_{1}A(U)+\xi_{2}B(U),
\end{equation}
where $\xi:=(\xi_{1},\xi_{2})^{T}$ is a unit vector. By calculations, the eigenvalues of $M(U)$ read
\begin{equation}
\lambda_{1}\equiv\cdots\equiv\lambda_{n}\equiv\lambda_{-},\;
\lambda_{n+1}\equiv\cdots\equiv\lambda_{2n}\equiv\lambda_{+},\;
\lambda_{2n+1}\equiv\cdots\equiv\lambda_{3n}\equiv0,
\end{equation}
where
\begin{equation}
\lambda_{\pm}=\frac{b\xi_{1}+c\xi_{2}\pm\sqrt{\Xi}}{a},
\end{equation}
in which
\begin{equation}
\Xi:=(b^{2}-ae)\xi_{1}^{2}+2(bc-ad)\xi_{1}\xi_{2}+(c^{2}-ap)\xi_{2}^{2}.
\end{equation}

By a direct calculation, the right eigenvector corresponding to $\lambda_{i}\,(i=1,\cdots,3n)$ can be chosen as
\begin{equation}
r_{i}=(-\lambda_{-}e_{i},\xi_{1}e_{i},\xi_{2}e_{i})^{T}\quad (i=1,\cdots,n),
\end{equation}
\begin{equation}
r_{i}=(-\lambda_{+}e_{i},\xi_{1}e_{i},\xi_{2}e_{i})^{T}\quad (i=n+1,\cdots,2n),
\end{equation}
\begin{equation}
r_{i}=(0,\xi_{2}pe_{i},-(e\xi_{1}+2\xi_{2}d)e_{i})^{T}\quad (i=2n+1,\cdots,3n),
\end{equation}
where, $e_{k}=(0,\cdots,0,\displaystyle{\mathop{1}^{(k)}},0,\cdots,0)\,(k=1,\cdots,n)$.

By careful calculations, $\lambda_{i}\,(i=1,\cdots,2n)$ are completely linearly degenerate in the sense of Definition 2.1 and $\lambda_{i}\,(i=2n+1,\cdots,3n)$ are linearly degenerate in the sense of P. D. Lax. For detailed computations, one can refer to \cite{HK}
\subsection{Lax system}
In this section, we give an example (see \cite{cg}) which satisfies Definition 2.1.

Let $f(u)$ be an analytic function of a single complex variable $U=u+iv$. We impose on the complex-valued function $U=U(t,z)$, $z=x+yi$ and the real variable $t$ of the following nonlinear PDE:
\begin{equation}
\partial_{t}\overline{U}+\partial_{z}f(U)=0,
\end{equation}
where the bar denotes the complex conjugate and $\partial_{z}=\frac{1}{2}(\partial_{x}-i\partial_{y})$. We may express this equation in terms of the real and imaginary parts of $U$ and $\frac{1}{2}f(U)=a(u,v)+b(u,v)i$. Then (4.19) gives
\begin{equation}
\left\{\begin{array}{l}\partial_{t}u+\partial_{x}a(u,v)+\partial_{y}b(u,v)=0,\vspace{2mm}\\
\partial_{t}v-\partial_{x}b(u,v)+\partial_{y}a(u,v)=0.\vspace{2mm}\\
\end{array}\right.
\end{equation}
By (4.20), Lax system is a hyperbolic conservation law. Thus, if system (4.20) is completely linearly degenerate in the sense of Definition 2.1, the solution of the Cauchy problem (4.19) with arbitrarily initial data exists globally.
\begin{Remark}
In our paper \cite{KW}, we have proved that under the assumptions of linearly degenerate in the sense of P. D. Lax, the life span of (4.19) depends on the flux function $f(U)$. In another paper \cite{WCH}, we derived a sharp blowup result under the assumptions that the system is genuinely nonlinear in some direction in the sense of P. D. Lax. While,
in the sense of Definition 2.1, we give a sufficient condition to guarantee the global existence of the Cauchy problem of (4.20).
\end{Remark}
\subsection{Quasilinear hyperbolic system of conservation laws with rotational
invariance}
    Consider the following quasilinear hyperbolic system
\begin{equation}
\frac{\partial u}{\partial t}+\sum_{i=1}^m\frac{\partial}{\partial x_i}
\left(f_i(|u|)u\right)=0,
\end{equation}
with initial data
$$
t=0:\qquad u(0,x)=u_{0}(x),
$$
where $u=(u_1,\dots,u_n)^T$, $x=(x_{1},\cdots,x_{m})$, $f_i\in C^2\left({\mathbb R}^{+},{\mathbb R}\right)$ and $m$ is an integer $\geq 1$. System (4.21) can be used to describe the
propagation of waves in various situations in mechanics (such as the reactive
flows, magneto-hydrodynamics and elasticity theory, etc.) at least for
the case that $m=1$ (\cite{B}, \cite{F1}, \cite{KK}, \cite{LW}).  It is no longer strictly
hyperbolic, and possesses the eigenvalues with constant multiplicity even for
the case that $m=1$.  When $m=1$ and $n=2$, system (4.21) was first studied
by \cite{KK} and \cite{LW}.  Freist\"{u}hler \cite{F1}-\cite{F2} considered the Riemann problem
and the Cauchy problem for system (4.21) with $m=1$ and  $n \ge 1$. When $m>1$ and $n>1$, system (4.21) was studied by Bressan \cite{AB}

Rewrite (4.21) as
\begin{equation}
\frac{\partial u}{\partial t}+\sum_{i=1}^m A_i(u)\frac{\partial u}
{\partial x_i}=0,
\end{equation}
where
\begin{equation}
A_i(u)=\left(\begin{array}{cccc}
f_i(r)+\frac{f_i'(r)}{r}u_1^2 & \frac{f_i'(r)}{r}u_1u_2
& \cdots & \frac{f_i'(r)}{r}u_1u_n \\[3pt]
\frac{f_i'(r)}{r}u_1u_2 & f_i(r)+\frac{f_i'(r)}{r}u_2^2
& \cdots & \frac{f_i'(r)}{r}u_2u_n \\[2pt]
\vdots & \vdots & \ddots & \vdots \\[2pt]
\frac{f_i'(r)}{r}u_1u_n & \frac{f_i'(r)}{r}u_2u_n &
\cdots & f_i(r)+\frac{f_i'(r)}{r}u_n^2
\end{array}\right),
\end{equation}
in which  $r=|u|>0$.

    In what follows, we consider the case that $r>0$.
Let
\begin{equation}
u=rs,
\end{equation}
where $r=|u|$, $s=(s_1,\cdots,s_n)^T\in S^{n-1}$. Then system (4.21) can be
rewritten as
\begin{equation}
\frac{\partial s}{\partial t}+\sum_{i=1}^m f_i(r)\frac{\partial s}
{\partial x_i}=0,\quad\;
\end{equation}
\begin{equation}
\frac{\partial r}{\partial t}+\sum_{i=1}^m\frac{\partial }{\partial x_i}
\left(rf_i(r)\right)=0.
\end{equation}
It is easy to check that system (4.21) is equivalent to system (4.25)-(4.26) at least for classical solutions. Moreover,
$$A(u,\xi)=\left(\begin{array}{cccccc}
\sum_{j=1}^mf_j(r)\xi_j & 0 & 0 & \cdots & 0 & 0 \\[3pt]
0 & \sum_{j=1}^mf_j(r)\xi_j & 0
& \cdots & 0 & 0 \\[2pt]
\vdots & \vdots & \vdots & \ddots & \vdots & \vdots\\[2pt]
0 & 0 & 0 & \cdots &\sum_{j=1}^mf_j(r)\xi_j  & 0 \\[2pt]
0 & 0 & 0 & \cdots & 0  & \sum_{j=1}^m(rf_j(r))^{\prime}\xi_j
\end{array}\right).$$
Obviously,
$\lambda(u,\xi)\stackrel{\triangle}{=}\sum_{j=1}^mf_j(r)\xi_j$ is an
eigenvalue of $A(u,\xi)$ with constant multiplicity $n$ be completely linearly
degenerate in the sense of Definition 2.1; while
$\mu(u,\xi)\stackrel{\triangle}{=}\sum_{j=1}^m(rf_j(r))^{\prime}\xi_j$
is genuinely nonlinear provided that for all $r>0$,
$\sum_{j=1}^m|(rf_j(r))^{\prime\prime}|>0$.

Consider the Cauchy problem (4.25) and (4.26) with the following initial data
\begin{equation}
t=0:\qquad r=|u_{0}(x)|,\quad s=\frac{u_{0}(x)}{|u_{0}(x)|}.
\end{equation}
Then the following theorem holds obviously
\begin{Theorem}
If the scalar conservation law (4.26) with initial data (4.27) admits a unique global $\mathbb C^{1}$ solution $r=r(t,x)$ on $t\geq0$, then the solution of the Cauchy problem of (4.21) admits a unique global $\mathbb C^{1}$ solution.
\end{Theorem}

\vskip 4mm

\noindent{\Large {\bf Acknowledgements.}} This work was supported in part by the NNSF of China (Grant No.: 11271323), Zhejiang Provincial Natural Science Foundation of China (Grant No: Z13A010002) and a National Science and Technology Project during the twelfth five-year plan of China (2012BAI10B04).


\begin{thebibliography}{99}
\bibitem{B}  M. Ben-Artzi, {\it The generalized Riemann problem for
reactive flows}, J. Comput. Phys., {\bf 31}, 70-101 (1989).
\bibitem{bre} A. Bressan, {\it Contractive metrics for nonlinear hyperbolic systems}, Indiana Univ. Math. J., {\bf 37}, 409-421 (1988).
\bibitem{AB} A. Bressan, {\it Some remarks on multidimensional systems of conservation laws}, Rend. Mat. Acc. Linceis., {\bf 15}, 225-233 (2004)
\bibitem{chr} D. Christodoulou, {\it Global solutions of nonlinear hyperbolic equations for small initial data}, Comm. Pure Appl. Math., {\bf 39}, 267-282 (1986).
\bibitem{cg} G. Q. Chen, {\it Multidimensional conservation laws: overview, problems, and perspective}, Nonlinear conservation laws and applications, {\bf 153}, 23-72 (2011).
\bibitem{dk} H. H. Dai and D. X. Kong, {\it Nonlinear degree and partial stability for quailinear hyperbolic systems and the application to plane elastic waves in hypereleastic materials}, Physics Letters A., {\bf 289}, 313-322 (2001).
\bibitem{F1} H. Freist\"{u}hler, {\it Rotational degeneracy of hyperbolic
systems of conservation laws}, Arch. Rat. Mech. Anal., {\bf 113}, 39-64 (1991).

\bibitem{F2} H. Freist\"{u}hler, {\it On the Cauchy problem for a class
of hyperbolic systems of conservation laws}, J. Diff. Equs., {\bf 112}, 170-178 (1994).
\bibitem{HK} S. J. Huang and D. X. Kong, {\it Equations for the motion of relativistic torus in the Minkowski space $\mathbb R^{1+n}$}, J. Math. Phys., {\bf 48}, 083510 (2007)
\bibitem{KW} D. X. Kong and C. H. Wei, {\it Global existence and life span of smooth solutions to a class of complex conservation laws}, Nonlinear Analysis, {\bf 95}, 553-567 (2014).

\bibitem{WCH} D. X. Kong and C. H. Wei, {\it Blowup of smooth solutions of a class of complex conservation laws},J. Math. Anal. Appl., {\bf 406}, 464-474 (2013)

\bibitem{KK} B. L. Keyfitz and H. Kranzer, {\it A system of non-strictly
hyperbolic conservation laws arising in elasticity theory}, Arch. Rat. Mech.
Anal., {\bf 72}, 219-241 (1980).

\bibitem{lax} P. D. Lax, {\it Hyperbolic systems of conservation laws II},
Comm. Pure Appl. Math., {\bf 10}, 537-566 (1957).
\bibitem{lax1} P. D. Lax, {\it Hyperbolic systems of conservation laws in several space variables}, Current topics in partial differential equations, 327-341, Kinokuniya, Tokyo, 1986.
\bibitem{LW} T. P. Liu and J. Wang, {\it  On a non-strictly hyperbolic
system of conservation laws}, J. Diff. Equs., {\bf 57}, 1-14 (1985).
\bibitem{ma} A. Majda, Compressible Fluid Flow and System of
Conservation Laws in Several Space Variables,  Applied Mathematical Sciences {\bf 53},
Springer-Verlag, 1984.

\end{thebibliography}
\end{document}